\documentclass[preprint,12pt]{elsarticle}

\usepackage{graphicx,epsfig}
\usepackage{amssymb,amsmath,bm}
\usepackage{rotating}
\usepackage{float}
\floatstyle{ruled}
\newfloat{code}{tbh}{aux}[section]
\floatname{code}{Code}

\journal{Journal of Computational Science}

\begin{document}

\begin{frontmatter}



\title{Open-source tools for dynamical analysis of Liley's mean-field cortex model}


\author{Kevin R. Green}
\author{Lennaert van Veen}

\address{Faculty of Science, University of Ontario Institute of Technology, 2000 Simcoe Street North, Oshawa, L1H 7K4 Ontario, Canada}

\begin{abstract}
Mean-field models of the mammalian cortex treat this part of the brain as a two-dimensional excitable medium.
The electrical potentials, generated by the excitatory and inhibitory neuron populations, are described by nonlinear, coupled,
partial differential equations, that are known to generate complicated spatio-temporal behaviour. We focus on the
model by Liley {\sl et al.} (Network: Comput. Neural Syst. (2002) 13, 67-113). Several reductions of this model have
been studied in detail, but a direct analysis of its spatio-temporal dynamics has, to the best of our knowledge, never
been attempted before. Here, we describe the implementation of implicit time-stepping of the model and the tangent linear model, 
and solving for equilibria 
and time-periodic solutions, using the open-source library PETSc. By using domain decomposition for 
parallelization, and
iterative solving of linear problems, the code is capable of parsing some dynamics of a macroscopic slice of
cortical tissue with a sub-millimetre resolution.
\end{abstract}

\begin{keyword}
Mean-field modelling \sep hyperbolic partial differential equations \sep 
numerical partial differential equations \sep
35Q92 \sep 
65Y05 
\end{keyword}

\end{frontmatter}


\section{Introduction}
\label{sec:intro}
Models of cortical dynamics come in two main families: network models and mean-field models. The former describe many interacting
neurons, each with their own dynamical rules, while the latter describe electrical potentials, generated collectively
by many neurons, as continuous in space and time. These potentials can be thought of
as averages over a number of macrocolumns, groups of hundreds of thousands of neurons in columnar structures at the
surface of the cortex. The reason to abandon the description of individual
neurons and pass to the mean-field limit, in analogy to the thermodynamic limit in statistical physics, is twofold.
Firstly, the description and analysis of a substantial piece of the cortex with a network model is not tractable since it would contain
billions of neurons, and many times more connections between them. As demonstrated by recent publications,
such as by Izhikevich and Edelman \cite{izhi} or by the Blue Brain team \cite{bluebrain}, progress in super computing allows
for the simulation of ever larger neuronal networks, that reflect actual brain dynamics. However, it is hard
to see how the output of such models can be analysed, other than by purely statistical techniques.
In contrast, mean-field models can be analysed as infinite-dimensional dynamical systems. 

The second advantage of the mean-field approach is that the electrical potentials, which appears as dependent
variables, are directly related to the electroencephalograph (EEG) \cite{Nunez2006}. The EEG
is usually measured with electrodes on the scalp or, in exceptional circumstances, directly on the surface of
the brain. In either case, the measured signal is not that of individual neurons, but that of many neurons, spread out
over a few square centimetres or millimetres. Thus, the way the signals of
individual neurons are smeared out by the spatial averaging of mean-field modelling is similar to the way they are mixed up in EEG
measurements.
Because of the direct link between
the local mean potential and the EEG, mean-field models are sometimes called EEG models.

The origin of mean-field modelling lies in the nineteen seventies, when pioneers like Walter Freeman \cite{Freeman1975},
Wilson \& Cowan \cite{Wilson1972} and Lopes da Silva \cite{LopesdaSilva1974} started to model components
of the human cortex with continuous fields. Over the past four decades, mean-field models have been used
to study a range of open questions in neuroscience, such as the generation of the alpha rhythm, 
8-$13\,\mathrm{Hz}$ 
oscillations in the EEG (see, e.g., \cite{LopesdaSilva1974,Nunez1974}), epilepsy (see, e.g., \cite{breakspear,kramer,blenkinsop}) 
and anaesthesia \cite{Bojak2005}.
Also, they are used in models for sensorimotor control, pattern discrimination and target tracking \cite{quinton}.

Although mean-field models have been used in all these contexts, little analysis has been done on
their behaviour as spatially extended dynamical systems. In part, this is due to their staggering
complexity. The Liley model \cite{Liley2002} considered here, for instance, consists of fourteen
coupled Partial Differential Equations (PDEs) with strong nonlinearities, imposed by coupling between
the mean membrane potentials and the mean synaptic inputs. The model can be reduced to a system of
Ordinary Differential Equations (ODEs) by considering only spatially homogeneous solutions, and the
resulting system has been examined in detail using numerical bifurcation analysis
(see \cite{Frascoli2011} and references therein). In order to compute equilibria, periodic orbits and such
objects for the PDE model, we need a flexible, stable simulation code for the model and it linearization
that can run in parallel to scale up to a domain size of about $2500\,\mathrm{cm}^2$, the size of 
a full-grown human cortex. We also need efficient, iterative solvers for linear problems with
large, sparse matrices. In this paper, we will show that all this can be accomplished in the
open-source software package PETSc \cite{petsc-user-ref}. Our implementation consists of
a number of functions in C that will be available publicly \cite{public}.

The goal of this computational work is similar to that of Coombes {\sl et al.}, who analysed
``spots'': rotationally symmetric, localized solutions in a model of a single neuron population
in two dimensions \cite{coombes}. The study of such special solutions will help us parse the
spatio-temporal dynamics of mean-field models. We will attempt to compute periodic orbits
and other special solutions in a full-fledged, two-population mean field model without imposing
any spatial symmetries.

\subsection{Liley's model}\label{model}

The model we use was first proposed by Liley {\sl et al.}~in 2002 \cite{Liley2002}.  
The dependent variables are the mean inhibitory and excitatory membrane potential,
$h_i$ and $h_e$, the four mean synaptic inputs, originating from either population 
and connecting to either, $I_{ee}$, $I_{ei}$, $I_{ie}$ and $I_{ii}$, and the excitatory
axonal activity in long-range fibers, connecting to either population, $\phi_{ee}$ and $\phi_{ei}$.
The model equations are
\begin{equation} \label{eq:Liley1}
\tau_k \frac{\partial h_k(\vec{x},t)}{\partial t} = h^r_k - h_k(\vec{x},t) 
  + \frac{h^{eq}_{ek}-h_k(\vec{x},t)}{\left|h^{eq}_{ek}-h^r_e\right|}I_{ek}(\vec{x},t)
  + \frac{h^{eq}_{ik}-h_k(\vec{x},t)}{\left|h^{eq}_{ik}-h^r_e\right|}I_{ik}(\vec{x},t)
\end{equation}
\begin{equation} \label{eq:Liley2}
\left(\frac{\partial}{\partial t} + \gamma_{ek}\right)^2 I_{ek}(\vec{x},t) = 
   e \Gamma_{ek}\gamma_{ek}\left\{N^\beta_{ek}S_e\left[h_e(\vec{x},t)\right]
                               +p_{ek}+\phi_{ek}(\vec{x},t)\right\}
\end{equation}
\begin{equation} \label{eq:Liley3}
\left(\frac{d}{dt} + \gamma_{ik}\right)^2 I_{ik}(\vec{x},t) = 
   e \Gamma_{ek}\gamma_{ek}\left\{N^\beta_{ik}S_i\left[h_i(\vec{x},t)\right]
                               +p_{ik}\right\}
\end{equation}
\begin{equation} \label{eq:Liley4}
\left[\left(\frac{\partial}{\partial t}+v\Lambda\right)^2 
      -\frac{3}{2}v^2\nabla^2\right]\phi_{ek}(\vec{x},t)   = 
N^\alpha_{ek}v^2\Lambda^2 S_e\left[h_e(\vec{x},t)\right]
\end{equation}
\begin{equation} \label{eq:Liley5}
S_k\left[h_k\right] = S^{max}_k\left(
  1 + \exp\left[-\sqrt{2}\frac{h_k-\mu_k}{\sigma_k}\right]
           \right)^{-1}
\end{equation}
where index $k=\{e,i\}$ denotes excitatory or inhibitory.  The meaning of the  parameters, along with 
some physiological bounds and the values used in our tests, are given in Table \ref{table:pars}. A detailed description 
of these equations can be found in references \cite{Liley2002,Frascoli2011}. Here, we will focus
on the aspects of the model most relevant for the numerical implementation. 

There are two sources of nonlinearity, related to the coupling of the synaptic inputs to the membrane
potential and vice versa. The former connection is quadratically nonlinear, while the latter is given
by the sigmoidal function $S$, which describes the onset of firing as the potential exceeds the threshold
value $\mu_{i,e}$. These nonlinearities tend to form sharp transitions of the potentials across the domain.
That is one reason why we opted for a finite-difference discretization over a pseudo-spectral approach.
Spectral accuracy would be of limited value in the presence of steep gradients and the finite-difference
scheme can be parallelized much more efficiently. The second reason is that we would like to be
able to change the geometry of the domain and the boundary conditions in future work. The finite-difference
scheme is more flexible in that respect.

The only spatial derivatives in the model are those in the equations for the long-range connections.
These are damped wave equations. We will discretize the Laplacian using a five-point stencil on a rectangular
grid. In previous work, Bojak \& Liley chose a second-order centered difference scheme for the time 
derivatives \cite{Bojak2005}. A disadvantage of this approach is that the stability condition of this scheme dictates that we 
set the time step inversely proportional to the grid spacing. In practice, they used a time step of $0.05\,\mathrm{ms}$.
To avoid this obstacle, we implemented the unconditionally stable implicit Euler method, as described in Sec.~\ref{timestepping}.

Other authors have used this model with an additional diffusive term in the equations for the membrane
potentials to model gap junctions \cite{steynross}. Inclusion of these terms can drastically change
the bifurcation behaviour, as they can cause Turing transitions to space-dependent equilibria.
Without the additional terms, a Hopf bifurcation from a spatially homogeneous 
equilibrium to a space dependent periodic orbit or a saddle-node bifurcation of this equilibrium
can be the primary instability. The gap junction terms can readily be included in our implementation, and in Sec.~\ref{equilibria}
we will describe how to solve for equilibrium states that may depend on space.

\begin{sidewaystable}[htbf]
\begin{tabular}{llllll}

\hline\noalign{\smallskip}
Parameter & Definition & Minimum & Maximum & Value & Units  \\
\hline
$h^r_e$ & resting excitatory membrane potential & $-80$ & $-60$ & -72.293 & mV\\
$h^r_i$ & resting inhibitory membrane potential & $-80$ & $-60$ & -67.261 & mV\\
$\tau_e$ & passive excitatory membrane decay time & $5$ & $150$ & 32.209 & ms\\
$\tau_i$ & passive inhibitory membrane decay time & $5$ & $150$ & 92.260 & ms\\
$h^{\mathrm{eq}}_{ee}$ & excitatory reversal potential & $-20$ & $10$ & 7.2583 & mV\\
$h^{\mathrm{eq}}_{ei}$ & excitatory reversal potential & $-20$ & $10$ & 9.8357 & mV\\
$h^{\mathrm{eq}}_{ie}$ & inhibitory reversal potential & $-90$ & $ h^r_k-5$ & -80.697 & mV\\
$h^{\mathrm{eq}}_{ii}$ & inhibitory reversal potential & $-90$ & $ h^r_k-5$ & -76.674 & mV\\
$\Gamma_{ee}$ & EPSP peak amplitude & $0.1$ & $2.0$ & 0.29835 & mV\\
$\Gamma_{ei}$ & EPSP peak amplitude & $0.1$ & $2.0$ & 1.1465 & mV\\
$\Gamma_{ie}$ & IPSP peak amplitude & $0.1$ & $2.0$ & 1.2615 & mV\\
$\Gamma_{ii}$ & IPSP peak amplitude & $0.1$ & $2.0$ & 0.20143  &mV\\
$\gamma_{ee}$& EPSP characteristic rate constant$^\ddagger$ & $100$ & $1,000$ & 122.68 & $\mathrm{s}^{-1}$\\
$\gamma_{ei}$& EPSP characteristic rate constant$^\ddagger$ & $100$ & $1,000$ & 982.51 & $\mathrm{s}^{-1}$\\
$\gamma_{ie}$& IPSP characteristic rate constant$^\ddagger$ & $10$ & $500$ & 293.10 & $\mathrm{s}^{-1}$\\
$\gamma_{ii}$& IPSP characteristic rate constant$^\ddagger$ & $10$ & $500$ & 111.40 & $\mathrm{s}^{-1}$\\
$N^\alpha_{ee}$ & no.\ of cortico-cortical synapses, target excitatory & $2000$ & $5000$ & 3228.0 & --\\
$N^\alpha_{ei}$ & no.\ of cortico-cortical synapses, target inhibitory & $1000$ & $3000$ & 2956.9 & --\\
$N^\beta_{ee}$ & no.\ of excitatory intracortical synapses & $2000$ & $5000$ & 4202.4 &  --\\
$N^\beta_{ei}$ & no.\ of excitatory intracortical synapses & $2000$ & $5000$ & 3602.9  & --\\
$N^\beta_{ie}$ & no.\ of inhibitory intracortical synapses  & $100$ & $1000$ & 443.71  & --\\
$N^\beta_{ii}$ & no.\ of inhibitory intracortical synapses  & $100$ & $1000$ & 386.43  & --\\
$v$ & axonal conduction velocity & $100$ & $1,000$ & 116.12 & $\mathrm{cm}\,\mathrm{s}^{-1} $\\
$1/\Lambda$ & decay scale of cortico-cortical connectivity  & $1$ & $10$ & 1.6423 & cm \\
$S^{\mathrm{max}}_e$ & maximum excitatory firing rate & $50$ & $500$  & 66.433  & $\mathrm{s}^{-1}$ \\
$S^{\mathrm{max}}_i$ & maximum inhibitory firing rate & $50$ & $500$  & 393.29  & $\mathrm{s}^{-1}$ \\
$\mu_e$ & excitatory firing threshold &  $-55$ & $-40$ & -44.522 & mV \\
$\mu_i$ & inhibitory firing threshold &  $-55$ & $-40$ & -43.086 & mV \\
$\sigma_e$ & standard deviation of excitatory firing threshold  & $2$ & $7$ & 4.7068  & mV \\
$\sigma_i$ & standard deviation of inhibitory firing threshold  & $2$ & $7$ & 2.9644  & mV \\
$p_{ee}$ & extracortical synaptic input rate & $0$ & $10,000$ & 2250.6 & $\mathrm{s}^{-1}$ \\
$p_{ei}$ & extracortical synaptic input rate & $0$ & $10,000$ & 4363.4 & $\mathrm{s}^{-1}$ \\
\hline
\end{tabular}
\caption{{\bf Meaning, ranges and values for the model parameters.}
The values used for the tests presented in Sec. \ref{examples} are taken from reference \cite{Bojak2007}.}
\label{table:pars}
\end{sidewaystable}
We will test our implementation by comparing to, and extending, the computations of oscillations
with a $40\,\mathrm{Hz}$ 
component by Bojak \& Liley \cite{Bojak2007}. The corresponding parameter values
are listed in Table~\ref{table:pars}. The $40\,\mathrm{Hz}$ 
oscillations arise spontaneously if the number of
local inhibitory-to-inhibitory connections is changed slightly. We introduce a scaling parameter
$r$ by replacing $N^{\beta}_{ii}\rightarrow r N^{\beta}_{ii}$. This is the only parameter that will
be varied in our tests.

\subsection{PETSc overview}
Rather than creating our code from scratch, we opted to work with the 
Portable, Extensible Toolkit for Scientific Computation (PETSc): an open-source, object 
oriented library that is designed for the scalable
solution and analysis of PDEs \cite{petsc-web-page,petsc-user-ref}.  
PETSc is written in the C language, and is usable from C/C++ as 
well as Fortran and Python. We use PETSc in conjunction with the Scalable Library for 
Eigenvalue Problem Computations (SLEPc) \cite{Hernandez2005}, for the computation
of eigenspectra of equilibrium and periodic solutions. Since our implementation
uses some features of PETSc and SLEPc that are recent additions and are still being tested, 
we use the development version of both projects.

\begin{figure}[hbt]
\begin{center}
\epsfig{file=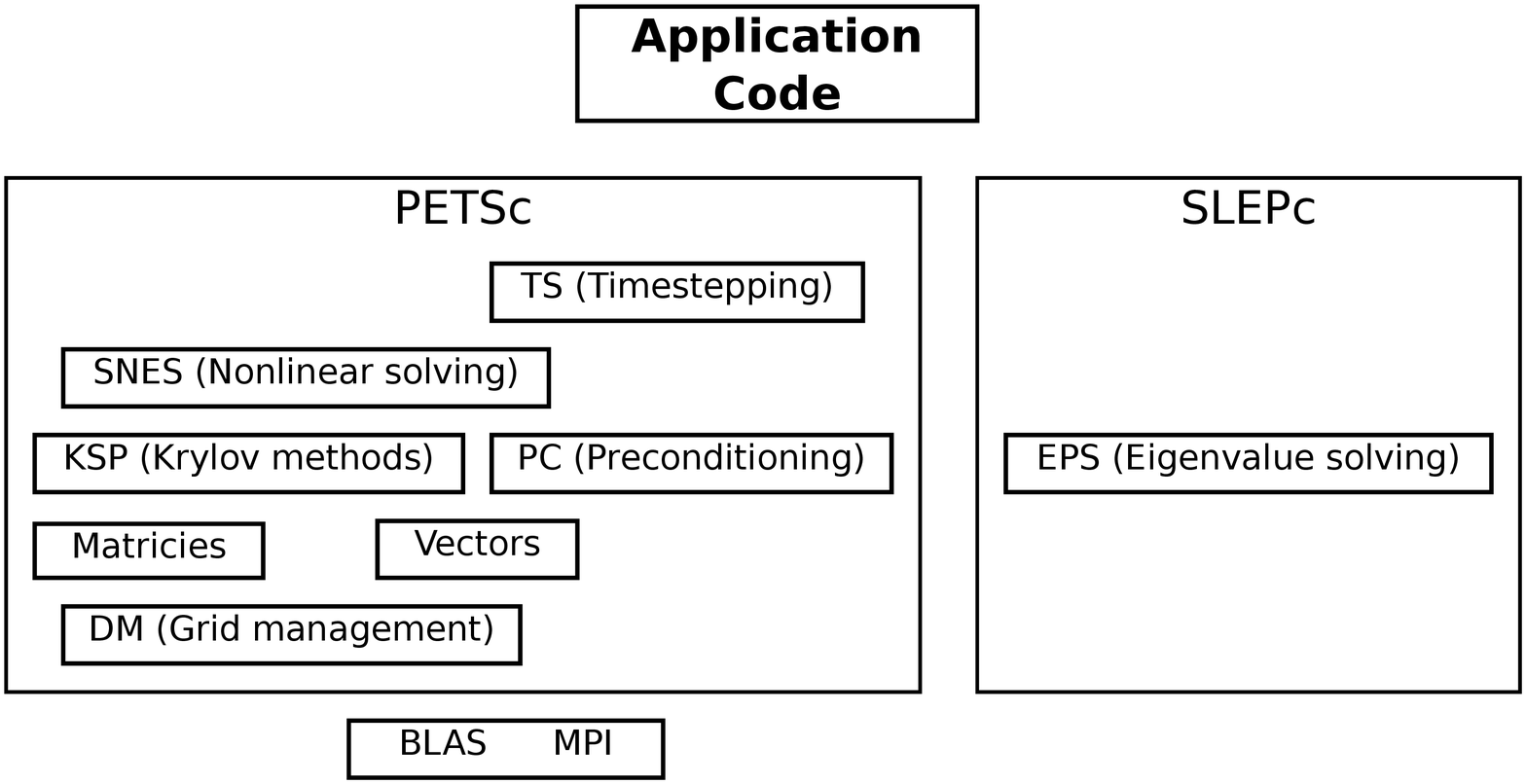,width=0.9\textwidth}
\caption{\label{fig:petsc-schematic} Schematic representation of the components of PETSc 
and SLEPc used in our code, and their relative hierarchy.}
\end{center}
\end{figure}

PETSc is split up into multiple components to address the various problems associated with 
solving PDEs numerically.  For our purposes, we treat the {\tt DM} component, which handles 
the topology of the discretization, as the most fundamental, from which we can easily derive 
memory allocation and communication for distributed vectors ({\tt Vec}) and matrices ({\tt Mat}).
With vectors and matrices, we can now solve linear systems, such as those that 
arise in Newton iteration
for implicit time-stepping and the computation of equilibria and periodic orbits. PETSc's component for this is 
called {\tt KSP}, and it has numerous iterative solvers implemented, as well 
as preconditioners, ({\tt PC}), to increase convergence rates.  For implicit time-stepping, for example,
we use GMRES , preconditioned with incomplete LU (ILU) factorization, combined
with the block Jacobi method \cite{Saad1986,Saad2003}.
On top of the linear solvers come the nonlinear solvers, PETSc's {\tt SNES} component, which 
implements a few different methods, such as globally convergent 
Newton iteration with line search \cite{schnabel}.
Finally, PETSc provides a timestepping component, {\tt TS}, to obtain time dependant 
solutions.  Implemented here are numerous explicit and implicit schemes such as adaptive 
stepsize Runge-Kutta and implicit Euler. The implicit schemes use the {\tt KSP} component. 
A schematic of the hierarchy discussed here can be found in Fig.~\ref{fig:petsc-schematic}.

For our dynamical systems calculations we will frequently need to compute specific 
eigenvalues and eigenvectors for system-sized matrices.  
For this end, we use SLEPc, which implements 
iterative eigenvalue solvers 
using PETSc {\tt Vec} and {\tt Mat} distributed data structures.  The component of 
SLEPc that we use is {\tt EPS}, which has a few algorithms for iteratively 
solving eigenproblems. Its default algorithm is Krylov-Schur 
iteration.

\section{Model Implementation}

\subsection{Geometry}
Following earlier work by Bojak \& Liley (e.g. \cite{Bojak2005,Bojak2007}) we consider
the PDEs on a rectangular domain with periodic boundary conditions. On this domain,
we use a rectangular grid of $N_x$ by $N_y$ points.
In the tests presented
in Sec.~\ref{examples}, the domain and the grid are square. PETSc allows for more complicated
domain shapes and grids, that can be encoded in the {\tt DM} component, independent of the
higher-level components. 

Within {\tt DM}, PETSc provides a simpler subcomponent, {\tt DMDA}, 
for working with finite differences on structured grids such as our rectangle.  
If we specify a stencil to use for the spatial derivatives in the {\tt DMDA},
PETSc will automatically handle numerous things for parallel execution, such as
memory allocation and the communication setup for distributed vectors and for 
the distributed Jacobian matrix.

\subsection{Fields}
To make use of PETSc's solvers, the model must be written as a system of equations that is 
first order in time. This we achieve by introducing new states $J_{jk}$ and $\psi_{ek}$ 
according to
\begin{equation}\label{eq:1storder_1}
\frac{\partial I_{jk}}{\partial t} = J_{jk}-\gamma_{jk}I_{jk}
\end{equation}
\begin{equation}\label{eq:1storder_2}
\frac{\partial J_{jk}}{\partial t}  = 
e \Gamma_{jk}\gamma_{jk}\left\{N^\beta_{jk}S_j\left[h_j\right] + \phi_{jk}
                               +p_{jk}\right\} - \gamma_{jk}J_{jk}
\end{equation}
\begin{equation}
\frac{\partial \phi_{ek}}{\partial t} = \psi_{ek}-v^2\Lambda^2\phi_{ek}
\end{equation}
\begin{equation}\label{eq:1storder_3}
\frac{\partial \psi_{ek}}{\partial t}  = v^2\Lambda^2N^\alpha_{ek}S_e\left[h_e\right] + 
\frac{3}{2}v^2\nabla^2\phi_{ek}- v^2\Lambda^2\psi_{ek},
\end{equation}
with indices $j,k=\{e,i\}$.

We opted to use a struct, seen in Code.~\ref{fig:fieldstruct}, to store the fields, 
rather than a triply indexed array.
\begin{code}[b]
\small
\begin{verbatim}
typedef struct _Field{
  PetscReal h_e, h_i, 
    I_ee, J_ee, 
    I_ie, J_ie,
    I_ei, J_ei, 
    I_ii, J_ii, 
    phi_ee, psi_ee, 
    phi_ei, psi_ei;
} Field;
\end{verbatim}\vspace{-10pt}
\caption{\label{fig:fieldstruct}Struct for the fields.}
\end{code}
This allows the code to be more readable in the function and Jacobian evaluation routines.  
For example, one accesses the $\phi_{ee}$ component at grid point $(x_i,y_j)$ simply as 
{\tt u[j][i].phi\_ee}, provided that the elements of the array ({\tt Field **u;}) are stored
on the processor in which the call is made.

\subsection{Parameters}
All of the model parameters are stored in a struct designated as the application context.  
The application context is how PETSc gets problem related parameters into the user-defined 
functions needed by its solvers.  
\begin{code}
\small
\begin{verbatim}
typedef struct _AppCtx{
  PassiveReal hr_e,  hr_i, 
    tau_e, tau_i, 
    heq_ee, heq_ie,
    heq_ei, heq_ii,
    Gamma_ee, Gamma_ie, 
    Gamma_ei, Gamma_ii,
    gamma_ee, gamma_ie,
    gamma_ei, gamma_ii,
    Nalpha_ee, Nalpha_ei,
    Nbeta_ee, Nbeta_ie,
    Nbeta_ei, Nbeta_ii,
    v, Lambda, 
    Smax_e, Smax_i,
    mu_e, mu_i,
    sigma_e, sigma_i,
    p_ee, p_ei, 
    p_ie, p_ii;
    ...
} AppCtx;
\end{verbatim}\vspace{-10pt}
\caption{\label{fig:appctxstruct}Application context struct with the model parameters.}
\end{code}
Similar to the fields, this allows readable code for the parameters.  For example, 
one accesses the $\Gamma_{ie}$ parameter as {\tt user->Gamma\_ie}, if {\tt user} is defined 
as the pointer {\tt AppCtx *user;}. 
How the parameters show up in our struct for the application 
context is shown in Code~\ref{fig:appctxstruct}.

\subsection{User supplied functions}
In addition to the structs listed above, we need to provide PETSc with (at least)
a C function that computes the vector field for a given state.  We call this function 
{\tt FormFunction}, and from this PETSc is capable of approximating the Jacobian 
with various finite difference methods.  However,
we also supply a C function to explicitly compute the Jacobian, named {\tt FormJacobian}, 
because this allows for more efficient calculations, especially when looking at stepping 
the variational equations in Sec.~\ref{sec:variational}.

\section{Timestepping}\label{timestepping}

We use the implicit Euler method to time-step the discretized equations. As mentioned
in Sec.~\ref{model}, this allows us to take larger time steps than feasible with
explicit methods. Since we are aiming to compute
periodic orbits, rather than to generate long time series, the first order accuracy of
the method is not an issue. Once a periodic orbit is computed, the time-step size can
be reduced to increase accuracy. Another option is to use Richardson extrapolation to
increase the order of accuracy, using the same nonlinear solving as described below.

\subsection{Mathematical basis}
We symbolically write the dynamical system as
\begin{equation}\label{eq:ds}
\dot{u}=f(u),\quad f:\mathbb{R}^{N}\rightarrow\mathbb{R}^N.
\end{equation}
where $N$ is the total number of unknowns after discretization, in our case $14\times N_x\times N_y$.
The implicit Euler scheme for time integration is given by
\begin{equation}\label{eq:beuler}
u_{n+1} = u_{n} + \mathrm{d}t\,f(u_{n+1})
\end{equation}
where the subscript represents the step number, $\mathrm{d}t$ the step size, and $u_{0}$ 
the initial conditions. This nonlinear equation is solved by Newton iteration:
\begin{equation}\label{eq:newton_beuler}
u^{k+1}_{n+1}= u^k_{n+1}+ \mathrm{d}u^k, 
\end{equation}
where the superscript denotes the Newton iterate, and 
$\mathrm{d}u^k$ is the solution to the linear system
\begin{equation}\label{eq:newton_beuler_lin}
\left(\mathbb{I}-\mathrm{d}t\left.\frac{\partial f}{\partial u}\right|_{u^k_{n+1}}\right)\mathrm{d}u^k
= \mathrm{d}t\,f(u^k_{n+1})-u^k_{n+1}+u^k_n,
\end{equation}
where $\partial f/\partial u$ denotes the $N\times N$ Jacobian matrix.
Provided that the initial approximation, $u^0_{n+1}$, is close enough 
to the actual solution of equation \eqref{eq:beuler}, this iteration should converge 
quadratically.  This is achieved by making the initial approximation the 
result of an explicit Euler step
\begin{equation}\label{eq:beuler_init}
u^0_{n+1}=u_n + \mathrm{d}t\,f(u_n).
\end{equation}

As we scale up the size of our problems, it becomes the linear solve in equation 
\eqref{eq:newton_beuler_lin} that takes most time.  This problem is handled by using 
GMRES to solve the linear system.  
For large time steps, the spectrum of the matrix in Eq.~\ref{eq:newton_beuler_lin} is spread out,
and we need to precondition it for iterative solving.  We make use ILU, which has shown 
to be reliable \cite{Sanchez2002,Saad1994} for this type of problem. If we use more than one processor, PETSc
uses distributed storage for the matrix, and combines ILU with block Jacobi preconditioning.

\subsection{Implementation}
PETSc provides a simple interface for timestepping in its {\tt TS} component.  The basic 
code required to set up a {\tt TS} is given in Code~\ref{fig:ts_setup}.  With a {\tt TS} 
set up like this, the timestepping parameters are set from command line arguments at run time.
For example, to do implicit Euler timestepping for $40.67\,\mathrm{ms}$ 
with a time step of $0.1\,\mathrm{ms}$, one needs to provide the arguments
\newline
{\tt -ts\_type beuler -ts\_dt 0.1 -ts\_final\_time 40.67}.
\newline
In this specific case, since the final time is not an integer number of timesteps, PETSc will 
step past it, and interpolate at the desired time.
\begin{code}[t]
\begin{verbatim}
  TS     ts;
  TSCreate(PETSC_COMM_WORLD,&ts);
  TSSetProblemType(ts,TS_NONLINEAR);
  TSSetExactFinalTime(ts);
  TSSetRHSFunction(ts,PETSC_NULL,FormFunction,&user);
  TSSetRHSJacobian(ts,J,J,FormJacobian,&user);
  TSSetFromOptions(ts);
  TSSolve(ts,u,PETSC_NULL)
\end{verbatim}\vspace{-10pt}
\caption{\label{fig:ts_setup}PETSc code for setting up and running the timestepping. 
{\tt FormFunction} 
and {\tt FormJacobian} are user provided functions that compute the rhs of equation \eqref{eq:ds}, 
and its Jacobian respectively.  {\tt J} is an appropriately allocated matrix to hold the Jacobian, 
and {\tt u} a vector to hold the solutions.}
\end{code}

\section{Stepping of the variational equations}\label{sec:variational}

\subsection{Mathematical basis}
The variational equations for the dynamical system are written as
\begin{equation}\label{eq:vareq}
\dot{v} = \left.\frac{\partial f}{\partial u}\right|_uv,\quad v\in \mathbb{R}^N
\end{equation}
and must be integrated simultaneously with the dynamical system \eqref{eq:ds}.
Solving the variational equations allow us to compute the stability of solutions, and is
also an essential ingredient for the treatment of boundary value problems such as 
those that arise in the computation of periodic orbits.

Performing implicit Euler timestepping on the variational equations \eqref{eq:vareq} 
requires solutions of the linear problems
\begin{equation}\label{eq:beuler_vareq}
\left(\mathbb{I}-\mathrm{d}t\left.\frac{\partial f}{\partial u}\right|_{u_{n+1}}\right)v_{n+1}
= v_n.
\end{equation}
Since we already have the Jacobian of the dynamical system at timestep $n+1$, stepping the 
variational equations requires only one additional $N\times N$ linear solve per time step.

\subsection{Implementation}
In PETSc, we implement the timestepping of the variational equations as a {\tt MATSHELL}.
A {\tt MATSHELL} allows users to define their own matrix type.  
Within a {\tt MATSHELL}, one needs to give a context for storing the relevant data and write 
 functions for the desired matrix operation.
 For example, we point the operation {\tt MATOP\_MULT} to 
a function that takes the initial state of the variational system $v(0)$ as input, 
and outputs the result $v(T)$ at the end of the timestepping.
The context we use for the time stepping of the variational equations is shown in 
Code \ref{fig:PeriodIntegrationCtx}.  The function we provide for {\tt MATOP\_MULT} 
works by first taking a step of the {\tt TS}, then loading the Jacobian computed 
from that step and solving equation \eqref{eq:beuler_vareq}.  This is repeated until 
the {\tt TS} reaches its end.
\begin{code}
\small
\begin{verbatim}
typedef struct _PeriodIntegrationCtx{
  // timestepping of the original eqn
  TS      ts;
  Mat     tsJac;
  Vec     initState,endState,fullSol;
  // additional requirements for variational eqn
  Mat     J,eye;
  KSP     ksp; 
} PeriodIntegrationCtx;
\end{verbatim}\vspace{-10pt}
\caption{\label{fig:PeriodIntegrationCtx}The {\tt MATSHELL} context for timestepping of the 
variational equations. The {\tt TS} holds the relevant info for stepping the dynamical 
system}
\end{code}

The {\tt MATSHELL} thus defined can be used by SLEPc for the iterative computation of eigen
pairs. In particular, we will use this approach to compute the Floquet multipliers of periodic 
orbits.

\section{Equilibria}\label{equilibria}
Having set up the function {\tt FormFunction} for the right hand side of the dynamical 
system, and its Jacobian computation {\tt FormJacobian}, also used for time integration, 
we can set up equilibrium calculations using PETSc's {\tt SNES} component with very little 
effort.
\subsection{Mathematical Basis}
Equilibrium solutions to the dynamical system \eqref{eq:ds} are solutions that satisfy
\begin{equation}\label{eq:equilibriumsolution}
f(u)=0.
\end{equation}
To solve this, we can set up a Newton iteration scheme
\begin{equation}\label{eq:newton_equilibrium}
u^{k+1} = u^{k}+\mathrm{d}u^k
\end{equation}
with $\mathrm{d}u$ coming from the solution of the linear system
\begin{equation}\label{eq:newton_equilibrium_lin}
\left.\frac{\partial f}{\partial u}\right|_{u^k}\mathrm{d}u^k = -f(u^k).
\end{equation}
As with the timestepping, if the initial guess is good enough this will converge 
quadratically provided that $\left.\frac{\partial f}{\partial u}\right|_{u^k}$ is 
nonsingular.  Unlike the case of time stepping, though, we do not always have a
way to produce an initial approximation that is good enough.  
For stable equilibrium solutions, we can 
use timestepping to get close to an equilibrium, but this will 
not work for unstable equilibria.  One possible solution is using 
globally convergent Newton methods. Using such methods we can find equilibria from
very coarse initial data, at the cost of computing many iterations. The line
search algorithm and the trust region approach (see, e.g. \cite{schnabel}) are implemented in the {\tt SNES} component.

Stability of equilibrium solutions follows from the spectrum of the Jacobian. Because of the spatial symmetries of
the model, these will mostly appear in groups. On a square domain, for instance, a single eigenvalue will be
associated with up to eight eigenvectors, with wavenumbers $(\pm k_x,\pm k_y)$ and $(\pm k_y,\pm k_x)$. 

\subsection{Implementation}
Setting up and using a nonlinear solver within PETSc is straightforward, as shown in Code~\ref{fig:snes_equilibrium}.  
The default algorithm used by {\tt SNES} is Newton's 
method with line search.

\begin{code}
\small
\begin{verbatim}
  SNES    snes;
  SNESCreate(PETSC_COMM_WORLD,&snes);
  SNESSetFunction(snes,r,FormFunctionSNES,&user);
  SNESSetJacobian(snes,J,J,FormJacobianSNES,&user);
  SNESSetFromOptions(snes);
  SNESSolve(snes,PETSC_NULL,u);
\end{verbatim}\vspace{-10pt}
\caption{\label{fig:snes_equilibrium}Code snippet for solving for equilibria. Vectors 
{\tt r} and {\tt u} are preallocated, with {\tt u} being the initial approximation, and 
{\tt J} a preallocated matrix for the Jacobian.}
\end{code}

\section{Periodic solutions}\label{sec:periodic}

The primary instability in  the  Liley model is often a Hopf bifurcation, and periodic orbits
have been shown to play an important role in the dynamics of ODE reductions of the model 
(e.g. \cite{Frascoli2011,VanVeen2006}). However, space dependent periodic orbits have not
previously been computed and studied. Using PETSc data structures for bordered matrices,
in conjunction with a {\tt MATSHELL} structure, we can solve for periodic orbits based on the 
time stepping described in Secs.~\ref{timestepping} and \ref{sec:variational}.
 
\subsection{Mathematical basis}
Periodic orbits solve the boundary value problem
\begin{equation}\label{eq:periodic}
F(u,T) = \phi(u,T)-u = 0,
\end{equation}
where $\phi$ is the flow of the dynamical 
system \eqref{eq:ds}, and $T$ is the period.
Our strategy for solving this equation is essentially that of Sanchez et al. \cite{Sanchez2004},
namely Newton iterations combined with unconditioned GMRES iteration.
Linearising Eq.~\ref{eq:periodic} gives
\begin{equation}\label{eq:periodic_diff}
\left(D_u\phi(u,T)-\mathbb{I}\right)\mathrm{d}u + f(\phi(u,T))\mathrm{d}T = -F(u,T),
\end{equation}
where $D_u\phi$ is a matrix of derivatives of the flow with respect to its initial condition.
Upon convergence, this is the monodromy matrix of the periodic orbit.
This results in $N$ equations in $N+1$ unknowns, which must be closed by a phase condition.
We opted for the use of a Poincar\'{e} plane involving one of the state variables, $u_k$:
\begin{equation}\label{eq:Poin_plane}
\phi_k(u,T)-C = 0,
\end{equation}
where $C$ is set appropriately, for instance to the time-mean value of $u_k$. This choice gives 
the following bordered system 
\begin{equation}\label{eq:periodic_bordered}
\left[
\begin{array}{cc}
\left(D_u\phi(u,T)-\mathbb{I}\right) & f(\phi(u,T)) \\
 \left[D_u\phi(u,T)\right]_{k,.} & f_k(\phi(u,T))
\end{array}\right]
\left[\begin{array}{c}
\mathrm{d}u \\
\mathrm{d}T
  \end{array}\right]
 = \left[\begin{array}{c} -F(u,T) \\ C-\phi_k(u,T)   \end{array}\right],
\end{equation}
where $[D_u\phi(u,T)]_{k,.}$ denotes the $k^{\rm th}$ row of the matrix $D_u\phi$.
An update can then be made to the approximate solution as 
\begin{equation}\label{eq:Newton_update}
\left[\begin{array}{c}
u^{n+1} \\
T^{n+1}
  \end{array}\right]
=\left[\begin{array}{c}
u^{n} \\
T^{n}
  \end{array}\right] + 
\left[\begin{array}{c}
\mathrm{d}u \\
\mathrm{d}T
  \end{array}\right].
\end{equation}

The matrix $D_u\phi$ is dense, so we should avoid calculating and storing it explicitly. 
Iterative solving of the linear problem, \eqref{eq:periodic_bordered}, requires the
computation of matrix-vector products, which are constructed from the integration of the 
variational equation \eqref{eq:vareq} with $v=\mathrm{d}u$ and the vector field $f(\phi(u,T))$
at the end point of the approximately periodic orbit.
Since the governing PDE is dissipative, most of the eigenvalues of the monodromy matrix
are clustered around zero. This aids the convergence of GMRES, without any preconditioning.
Sanchez et al.\ \cite{Sanchez2004} provide bounds for the number of GMRES iterations 
for the Navier-Stokes equation, and the convergence we observe for the Liley model is
qualitatively similar.

\subsection{Implementation}\label{implement:per}

The problem of creating a bordered matrix system in a distributed environment is 
not a trivial one.  The specific case that we have is one vector, $u$, 
that is sparsely connected and distributed among processors, and one parameter, 
$T$, that must exist and be synchronized across all processors.  

PETSc's {\tt DM} 
module has some recently introduced functionality that allows us to handle this in a 
straightforward way, letting us make use of the {\tt DMDA} already used in the 
other types of calculations.

{\tt DMRedundant} can be used for the $T$ component of our extended system, as it has 
the precise behaviour that we require.  Next, we use a {\tt DMComposite} to join together 
the {\tt DMDA} of the grid with the {\tt DMRedundant} of the period.  We can then derive 
vectors from this {\tt DMComposite}, and use these vectors for PETSc's iterative linear 
solvers.  PETSc code that illustrates this idea is shown in Code~\ref{fig:dmcomposite}.
\begin{code}
\begin{verbatim}
  DM  packer, redT;
  DMCompositeCreate(PETSC_COMM_WORLD,&packer);
  DMRedundantCreate(PETSC_COMM_WORLD,0,1,&redT);
  DMCompositeAddDM(packer,da);
  DMCompositeAddDM(packer,redT);
\end{verbatim}\vspace{-10pt}
\caption{\label{fig:dmcomposite}Additional {\tt DM} pieces for extended vectors as in equation 
\eqref{eq:Newton_update}, assuming that {\tt da} is the {\tt DM} associated with the 
grid structure.  The numerical arguments in {\tt DMRedundantCreate} represent the 
processor where the redundant entries live (in global vectors), and the number of 
redundant entries respectively.}
\end{code}

The matrix multiplication is done through a {\tt MATSHELL}, and the struct that holds the 
relevant data is found in Code~\ref{fig:periodrefine_struct}.

\begin{code}
\begin{verbatim}
typedef struct _PeriodFindCtx{
  Mat          *linTimeIntegration;
  DM           packer,redT;
  Vec          endState,f_at_endState;
} PeriodFindCtx;
\end{verbatim}\vspace{-10pt}
\caption{\label{fig:periodrefine_struct}For finding periodic solution, we need a method 
for integrating the variational equations (the {\tt MATSHELL} discussed in section 
\ref{sec:variational}), additional {\tt DM}s, 
and space for holding $f$ evaluated at the state at the end of the integration}
\end{code}

\section{Example calculations}
\label{examples}

In this section, we present some computations that serve to validate our implementation
and to investigate its efficiency. All tests are based on the parameter set in Table~\ref{table:pars},
and the scaling of the number of local inhibitory-to-inhibitory connections, $r$, is varied
around the first bifurcation from an equilibrium to more complicated, spatio-temporal behaviour.

\begin{figure}[hbt]
\begin{center}
\includegraphics[width=0.9\textwidth]{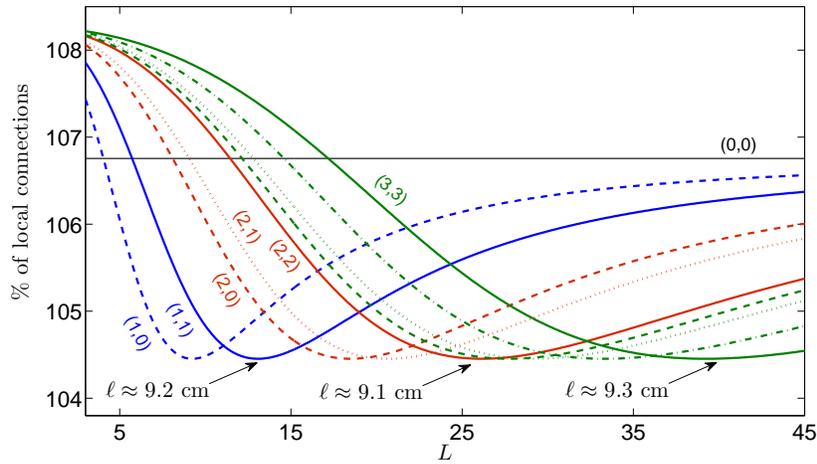}
\caption{\label{fig:nsc}Neutral stability curve for the spatially homogeneous equilibrium of the Liley model with
parameters set according to Table~\ref{table:pars}. Shown is the scaling parameter, $r$, versus the linear 
domain size, $L$, and wave numbers $\bm{k}=(k_x,k_y)$ are shown 
in parenthesis. When varying $r$, only for very small domains
the primary instability is spatially homogeneous. For domain sizes over $12.5\times 12.5\text{cm}^2$
the location of the primary instability approaches $r=1.04$ and the length scale of the leading instability
approaches $L/\|\bm{k}\|=9.3\,\mathrm{cm}$.}
\end{center}
\end{figure}

Fig.~\ref{fig:nsc} shows the neutral stability curve for the spatially homogeneous equilibrium,
which is the unique attractor of the model at small values of $r$. The primary transition is
a Hopf bifurcation with spatial wave numbers that depend on the system size. For systems smaller
than $2\times 2\,\mathrm{cm}^2$, the emerging periodic orbit is spatially homogeneous. For larger
systems, space dependent orbits emerge, and their typical length scale converges to about 
$9.3\,\mathrm{cm}$
for large system sizes. These stability curves were computed by solving small eigenvalue problems
for each combination of wavenumbers, independent from the PETSc implementation. The eigenvalues
computed by Krylov-Schur iteration in SLEPc, presented in Sec.~\ref{ex:eq}, are in good agreement.

\begin{figure}[hbt]
\epsfig{file=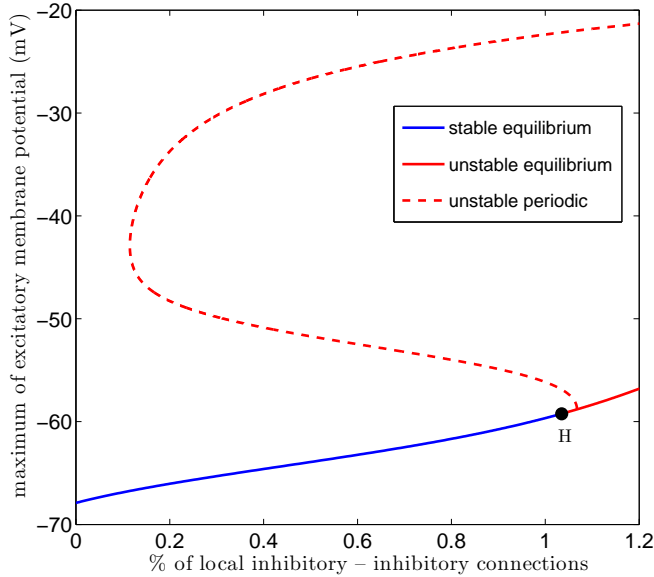,width=0.9\textwidth}
\caption{\label{fig:bifurcation_diagram}Partial bifurcation diagram showing the primary transition from
a spatially homogeneous equilibrium to a space and time dependent periodic orbit. On the vertical axis
the scaling parameter $r$ is plotted, and on the vertical axis the (maximum of) the excitatory membrane
potential. The branch of periodic solutions shown with a dashed line is a spatially homogeneous branch
that is unstable to space-dependent perturbations.}
\end{figure}

A partial bifurcation diagram, for spatially homogeneous solutions only, is shown in Fig.~\ref{fig:bifurcation_diagram}.
In this diagram, the Hopf bifurcation is subcritical, and time series analysis indicates that
the Hopf bifurcations associated with nonzero wave numbers are, too. The time series presented
in Sec.~\ref{ex:ts} was generated by starting from the equilibrium at $r=1$ and adding a finite-size
perturbation in the least stable direction, with wave number $|k_x|=|k_y|=1$.

\subsection{Timestepping}\label{ex:ts}

\begin{figure}[hbt]
\begin{center}
\includegraphics[width=.4\textwidth]{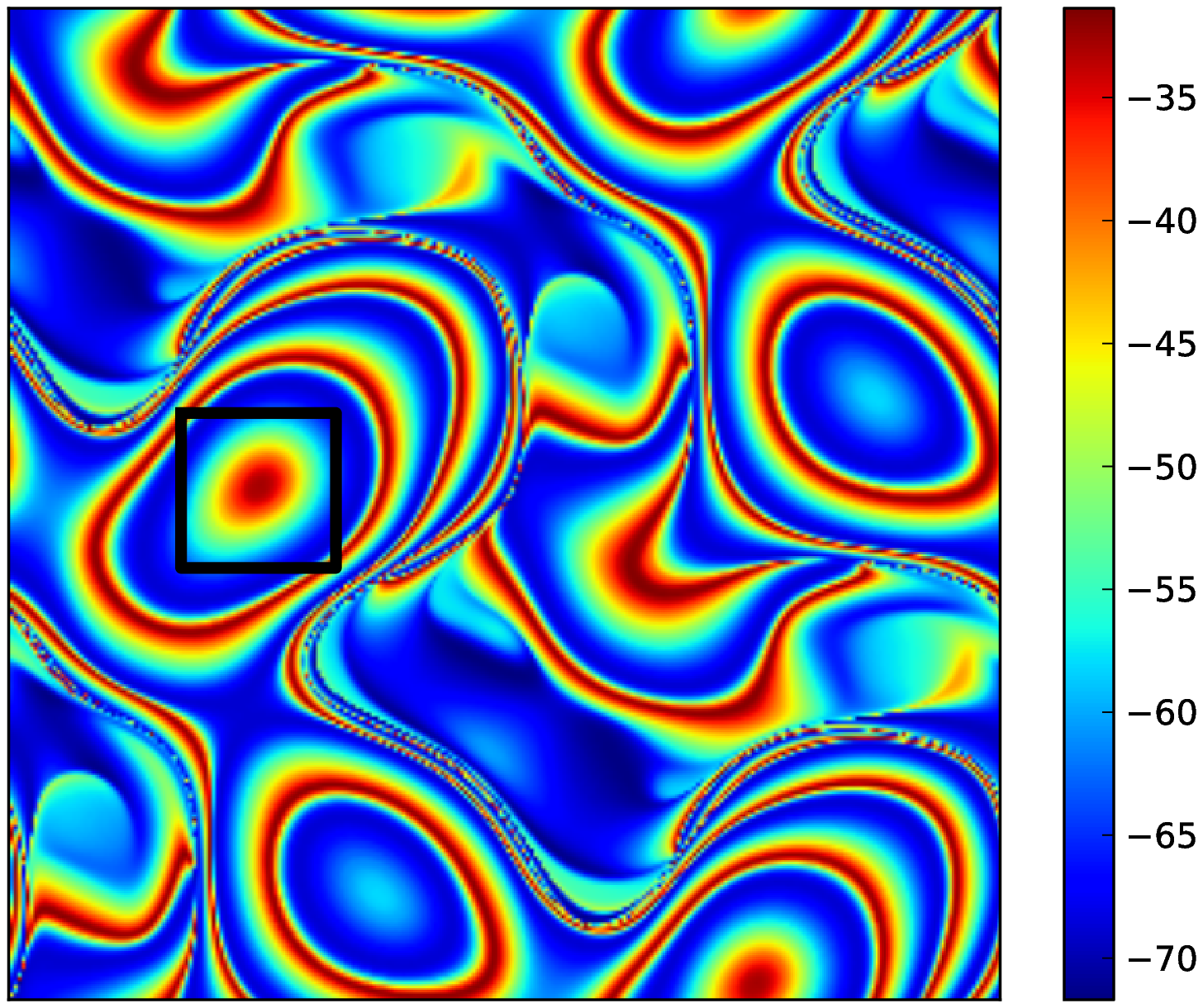}\qquad\quad\includegraphics[width=.4\textwidth]{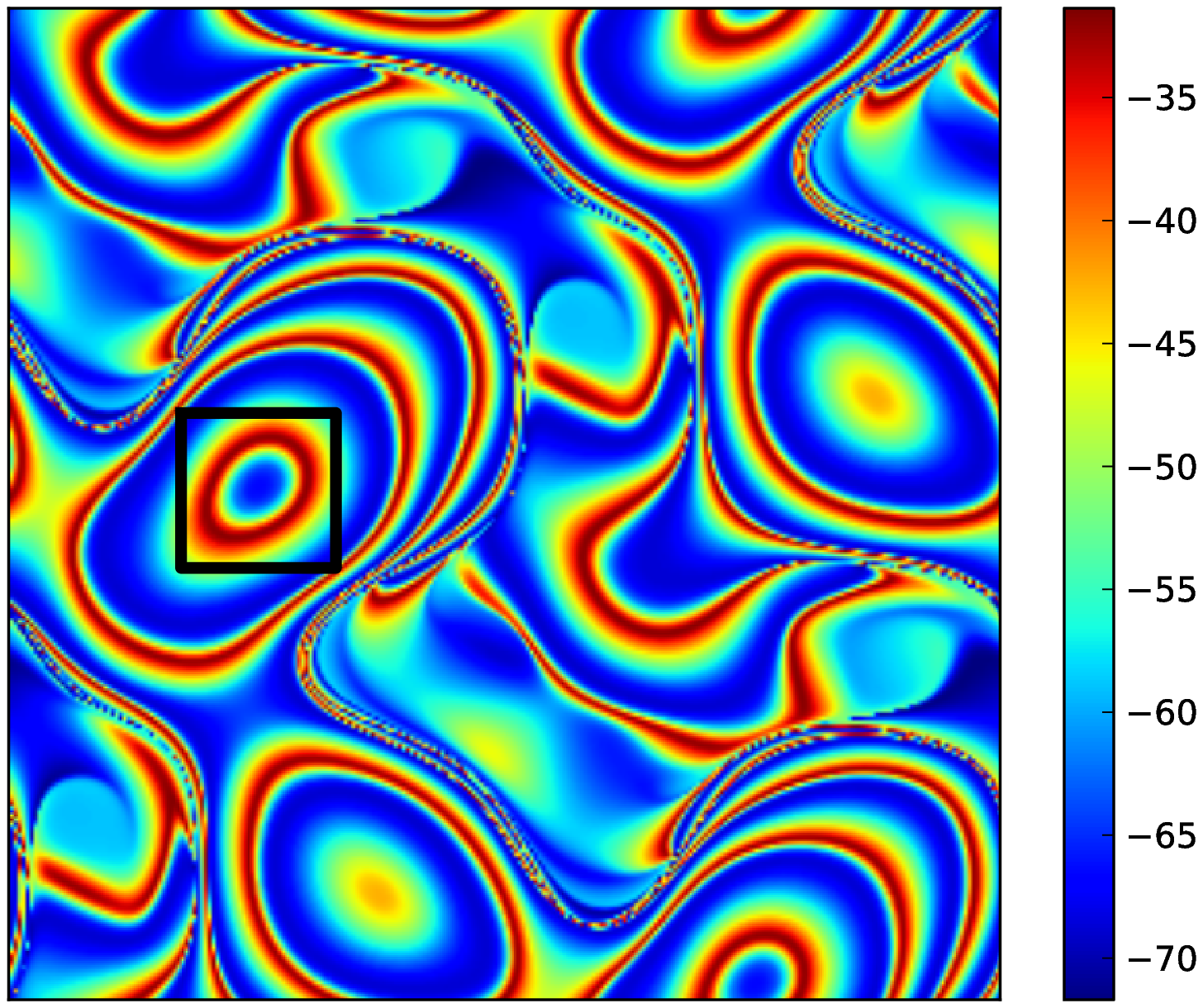}\vspace{10pt}\\
\includegraphics[width=.4\textwidth]{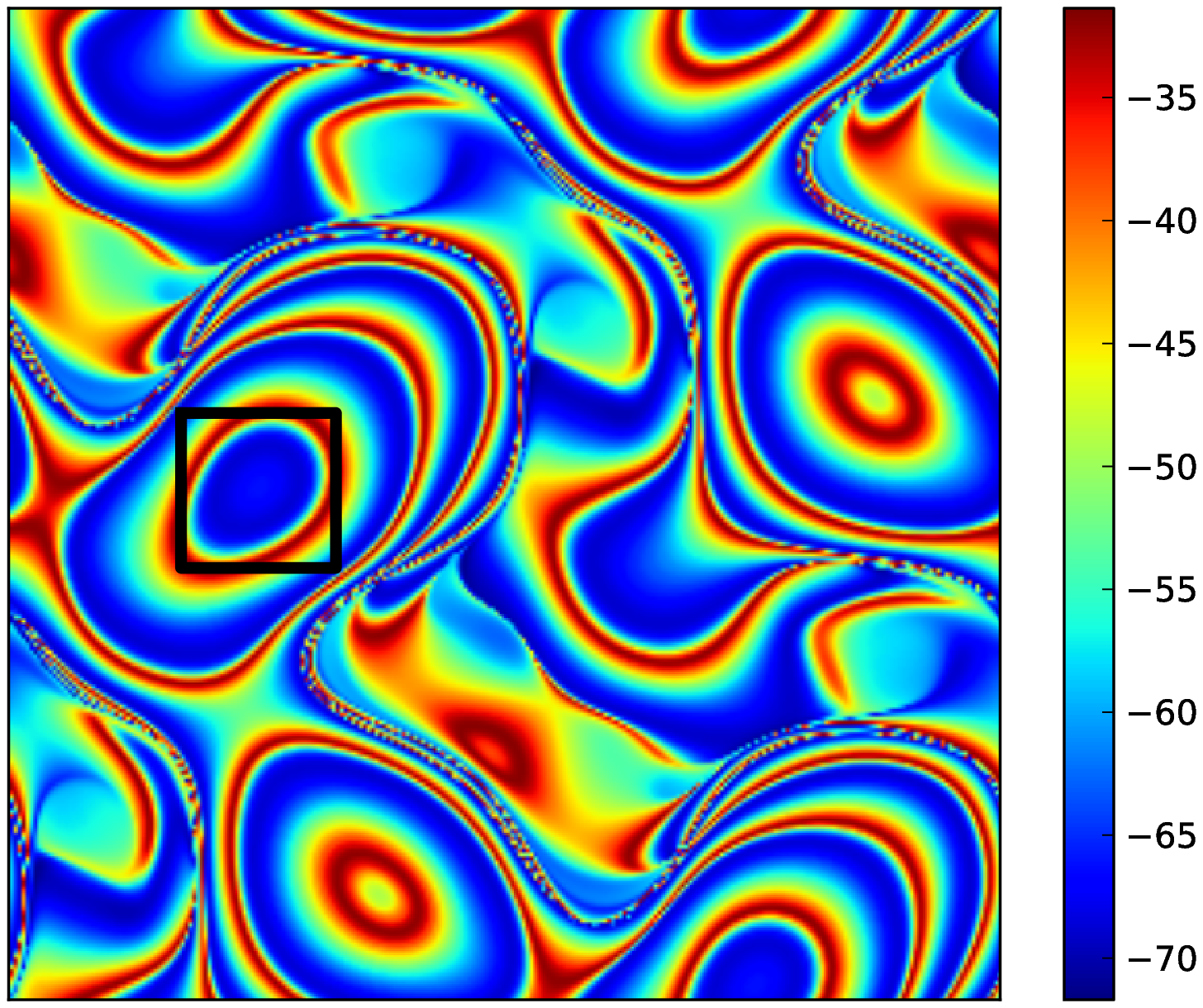}\qquad\quad\includegraphics[width=.4\textwidth]{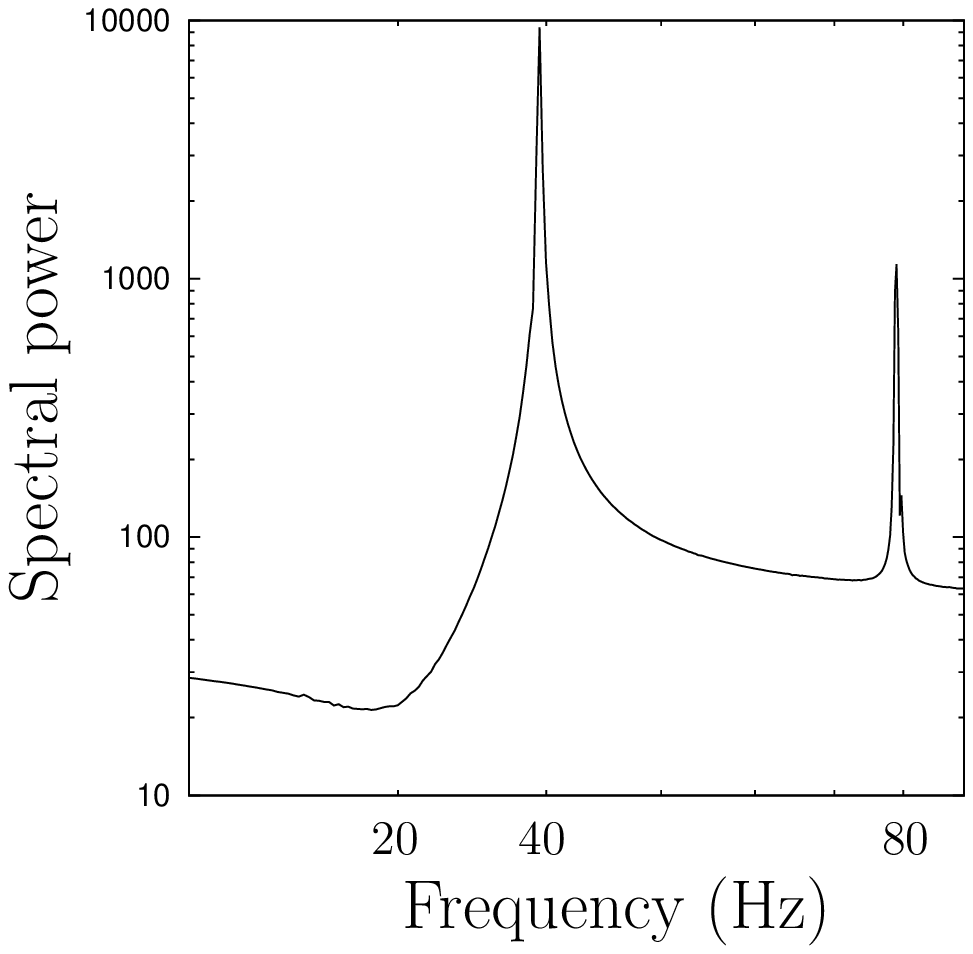}
\end{center}
\caption{\label{fig:time_snapshots}Three snapshots of the excitatory membrane potential, 
$6\,\mathrm{ms}$ apart, of a solution at $r=1$, near the primary Hopf 
bifurcation. The domain size is $12.8 \times 12.8\,\mathrm{cm}^2$, 
the resolution is $0.5\,\mathrm{mm}$ and the time-step 
size $1\,\mathrm{ms}$. 
The fourth panel shows the power spectrum of $h_e$, averaged over the region inside the black
square.}
\end{figure}

For the timestepping demonstration, we used a system size of $12.8\times 12.8\,\mathrm{cm}^2$ 
with $0.5\,\mathrm{mm}$ 
resolution, resulting in a $256\times 256$ grid, and $N=917,504$ unknowns in total.  Setting the 
parameter $r=1.0$, we initialize with the stable equilibrium solution perturbed by its least 
stable eigenmode, shown in Fig.~\ref{fig:eigenmode}.  Since the 
equilibrium solution is stable, small perturbations just decay, but sufficiently large perturbations
grow.  The snapshots 
of Fig.~\ref{fig:time_snapshots} were taken after a transient time of 600ms. The membrane potentials
show behaviour that is nearly periodic, with a dominant period of 40Hz, as demonstrated by the
power spectrum shown in the last panel.

Since the time-stepping code lies at the core of the periodic orbit solver, we carefully investigated its 
scaling with an increasing number of processors. Doubling the domain size, while keeping the grid spacing
fixed, gives a dynamical system with $N=3,670,016$ degrees of freedom. We time-stepped this system on 
a small subcluster of $2.4\,\mathrm{GHz}$ AMD Opteron nodes with gigabit interconnects. 
Apart from some minor load-balancing effects, the scaling is linear up to 16
processors, despite the relatively slow interconnects. 
\begin{figure}[hbt]
\epsfig{file=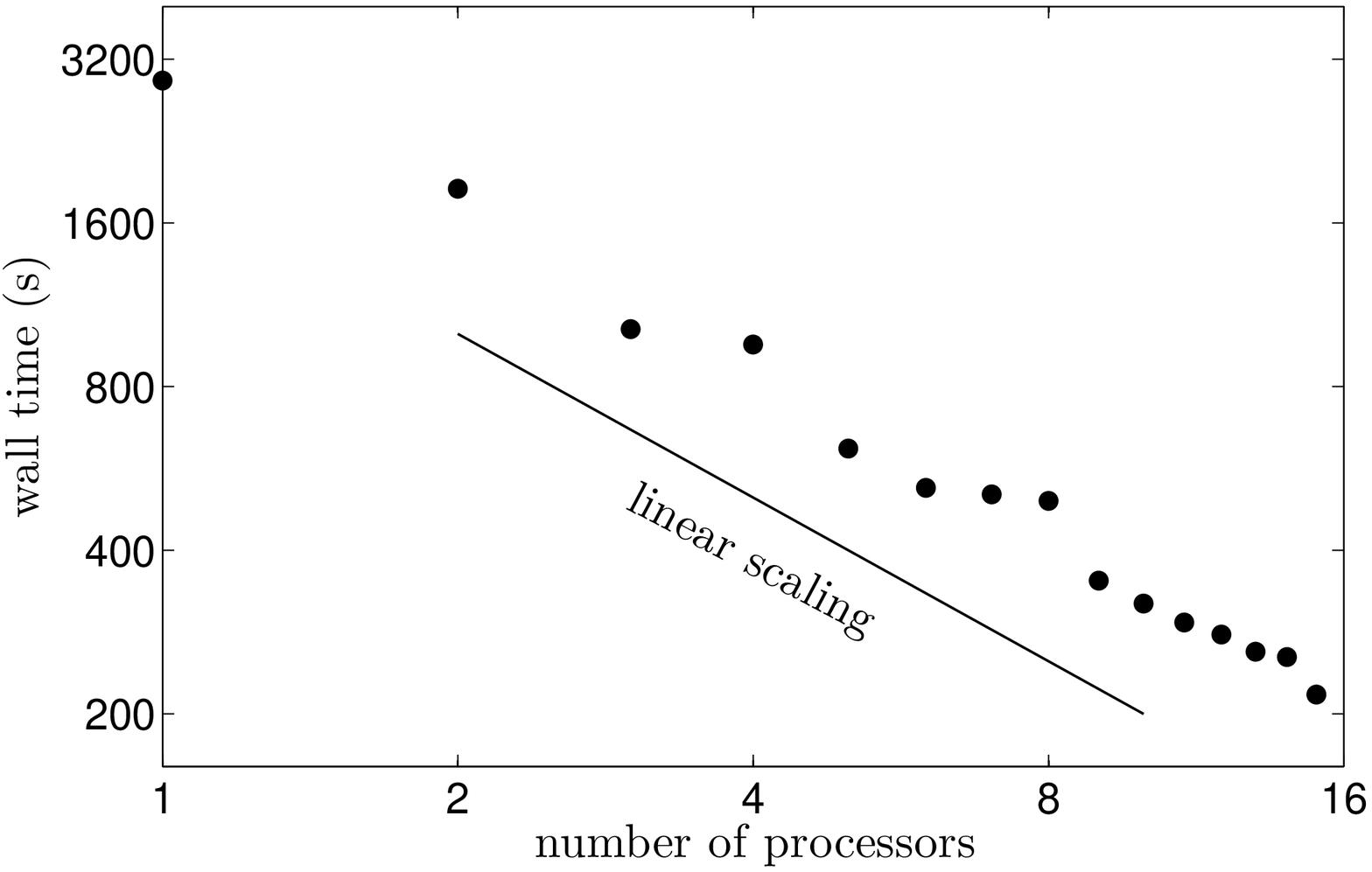,width=0.9\textwidth}
\caption{\label{fig:wtime}Wall time for the computation of 50 time steps of $1\,\mathrm{ms}$ each on a $25.6 \times 25.6\,\mathrm{cm}^2$
domain with $0.5\,\mathrm{mm}$ resolution. The fully implicit Euler steps are computed with Newton iterations, each of which
is solved for by GMRES, preconditioned with a combination of block Jacobi and ILU. The initial guess is given by 
an explicit Euler step. Two or three Newton iterations are sufficient to reduce the residual by a factor of $10^8$.
About 80 Krylov vectors are computed by GMRES to bring the relative residual down to $10^{-5}$. The number of unknowns 
is $N=3,670,016$.}
\end{figure}

\subsection{Equilibrium}\label{ex:eq}
\begin{figure}
\begin{center}
\includegraphics[width=.45\textwidth]{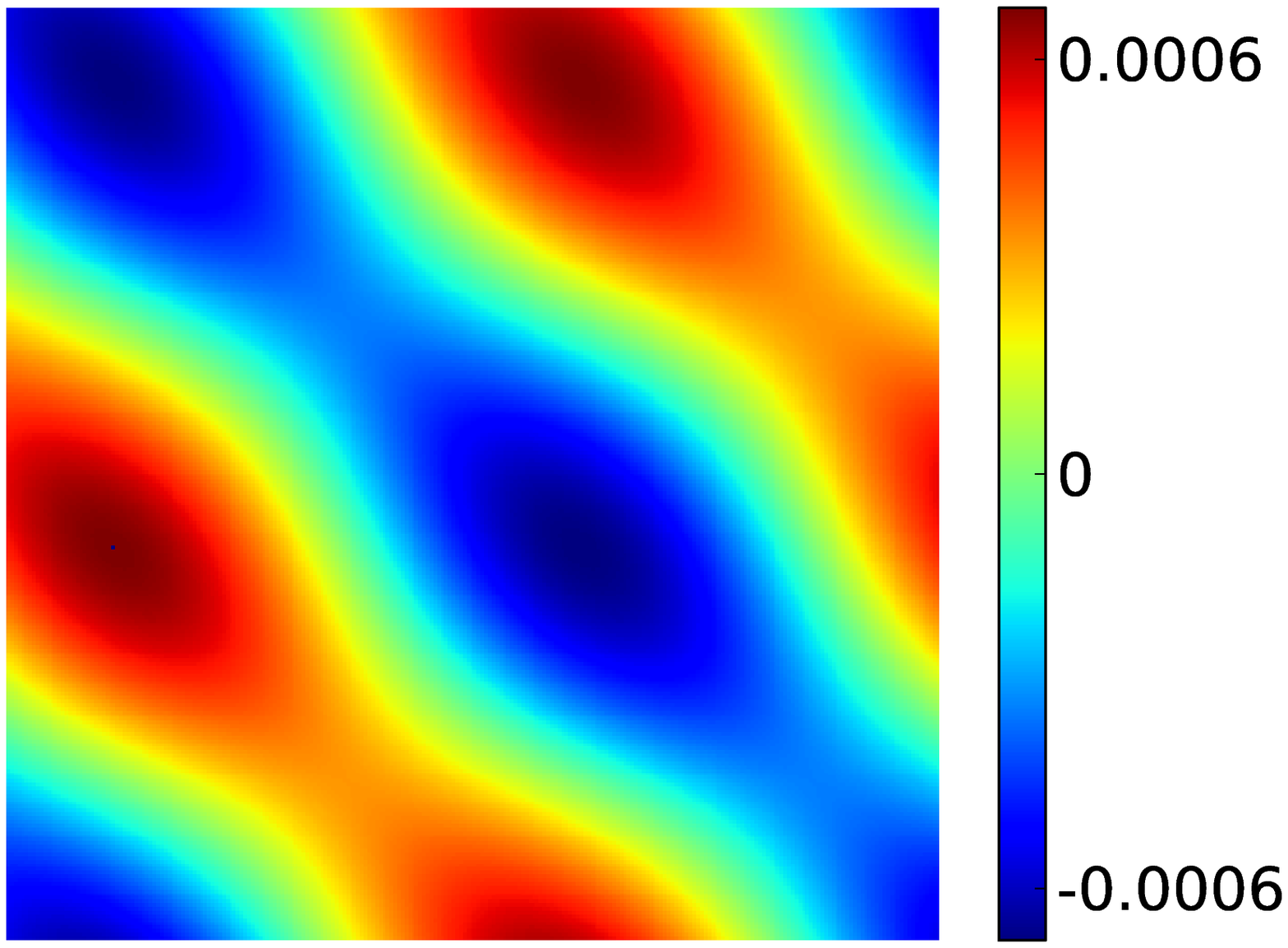}\qquad\includegraphics[width=.45\textwidth]{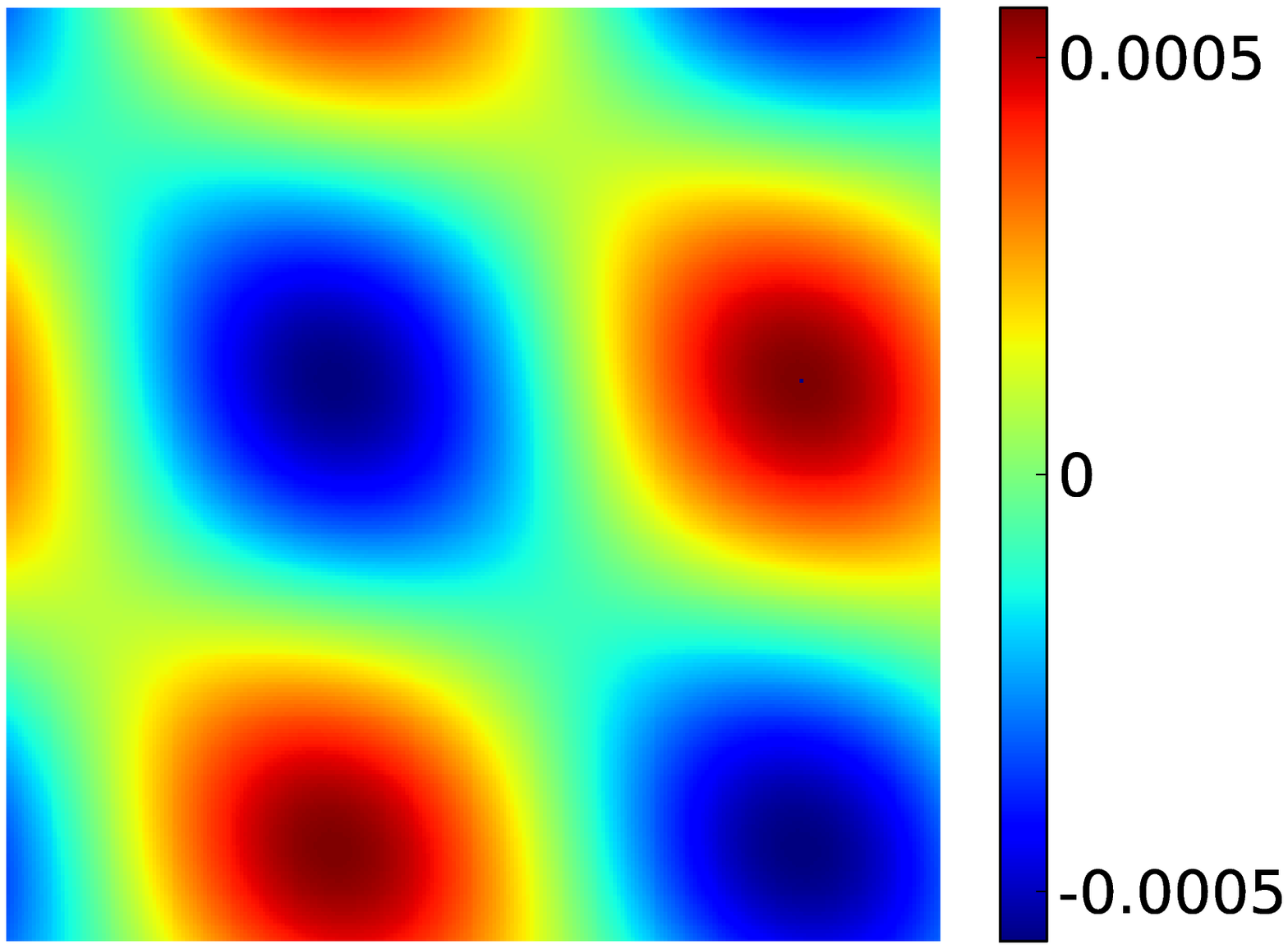}
\end{center}
\caption{\label{fig:eigenmode}The real part of the least stable eigenmode of the stable 
equilibrium located at $r=1.046$. Displayed are the excitatory (left) and inhibitory (right) membrane potentials.
The eigenvector, with wave numbers $(1,1)$, was computed by Arnoldi iteration and is scaled to have unit $L_2$ norm.}
\end{figure}

We computed the whole equilibrium curve of Fig.~\ref{fig:bifurcation_diagram} 
through parameter continuation, which is a trivial extension of the algorithm
for computing equilibria, presented in Sec.~\ref{equilibria}.  For each computed equilibrium 
solution, we took the Jacobian and used SLEPc to compute the 
eigenvalues with the largest real parts.  The result is shown in Fig.~\ref{fig:eigenmode}. As predicted
by the neutral stability curve computation, the $(1,1)$ mode turns unstable first, immediately followed by
the $(1,0)$ mode. Around $r=1.08$, the $(0,2)$ mode crosses the $(0,1)$ mode and proceeds to become the
most unstable mode for larger values  of $r$.
 The least stable 
eigenmode for $r=1.046$, just after its eigenvalue has crossed zero,
is shown in Fig.~\ref{fig:eigenmode}.
\begin{figure}[hbt]
\begin{center}
\includegraphics[width=0.8\textwidth]{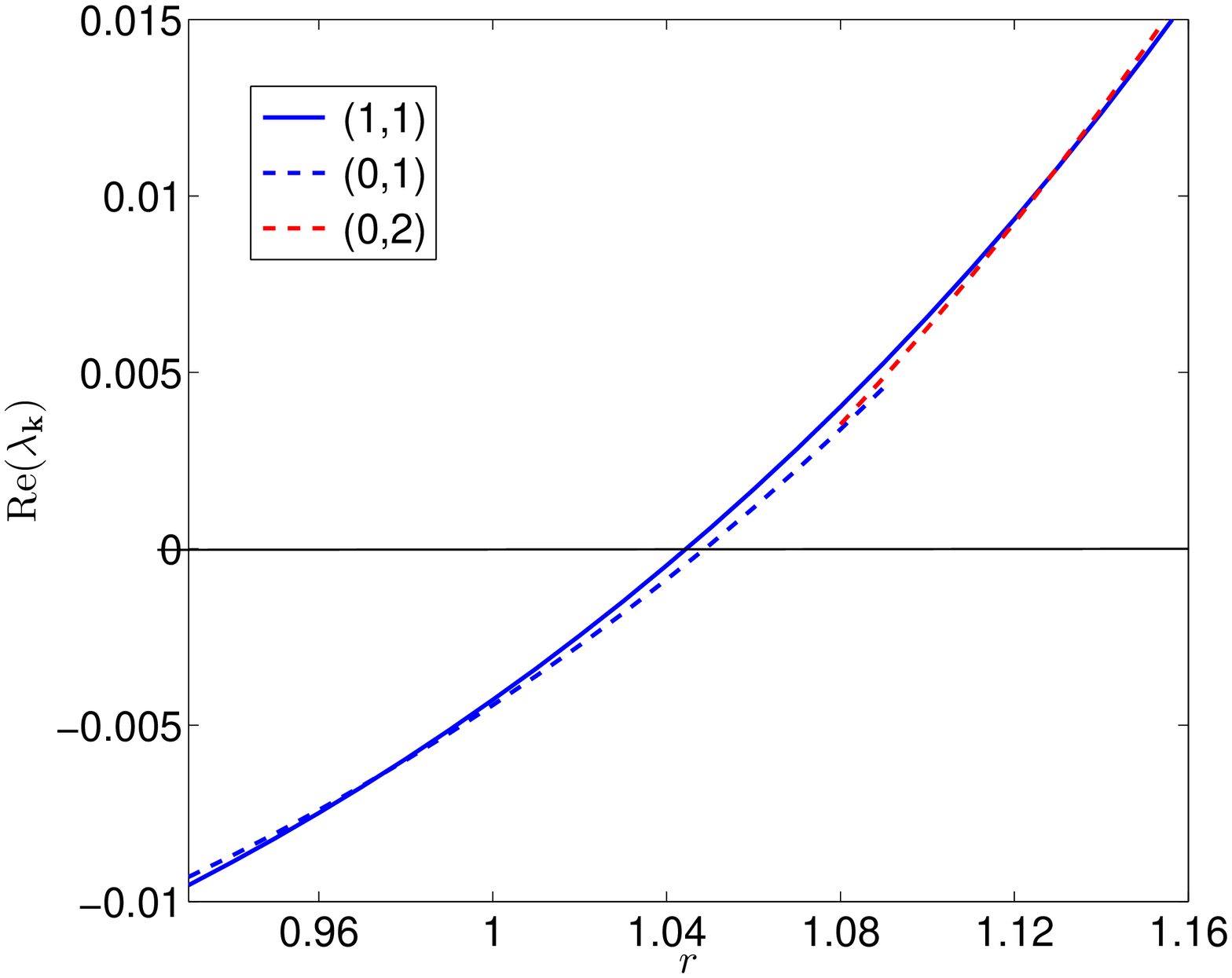}
\caption{\label{continuation}The real parts of the leading two eigenvalue pairs of the spatially homogeneous
equilibrium tracked in the scaling parameter $r$, for system size $L=12.8\,\mathrm{cm}$. The primary transition is tied 
to wave numbers $|k_x|=|k_y|=1$. The other curves shown are for wave numbers $k_x=0$, $k_y=\pm 1$ and
$k_x=\pm 1$, $k_y=0$ and for $k_x=0$, $k_y=\pm 2$ and
$k_x=\pm 2$, $k_y=0$.}
\end{center}
\end{figure}

\subsection{Periodic solutions}\label{ex:per}
We tested the computation of periodic orbits on a smaller grid, namely $16\times 16$ points, 
still with $0.5\,\mathrm{mm}$
resolution, and with $r=1.2$. The primary Hopf bifurcation is
sub critical, so there is no easy way to compute the branch of space-dependent periodic solutions.
Instead, we computed one of the spatially homogeneous orbits, for which an approximate solution 
can readily be obtained from analysis of the ODE reduction of the model. In fact, the upper part
of the branch of periodic orbits shown in Fig.~\ref{fig:bifurcation_diagram} is stable to all
spatially homogeneous perturbations.
\begin{figure}
\begin{center}
\includegraphics[width=.45\textwidth]{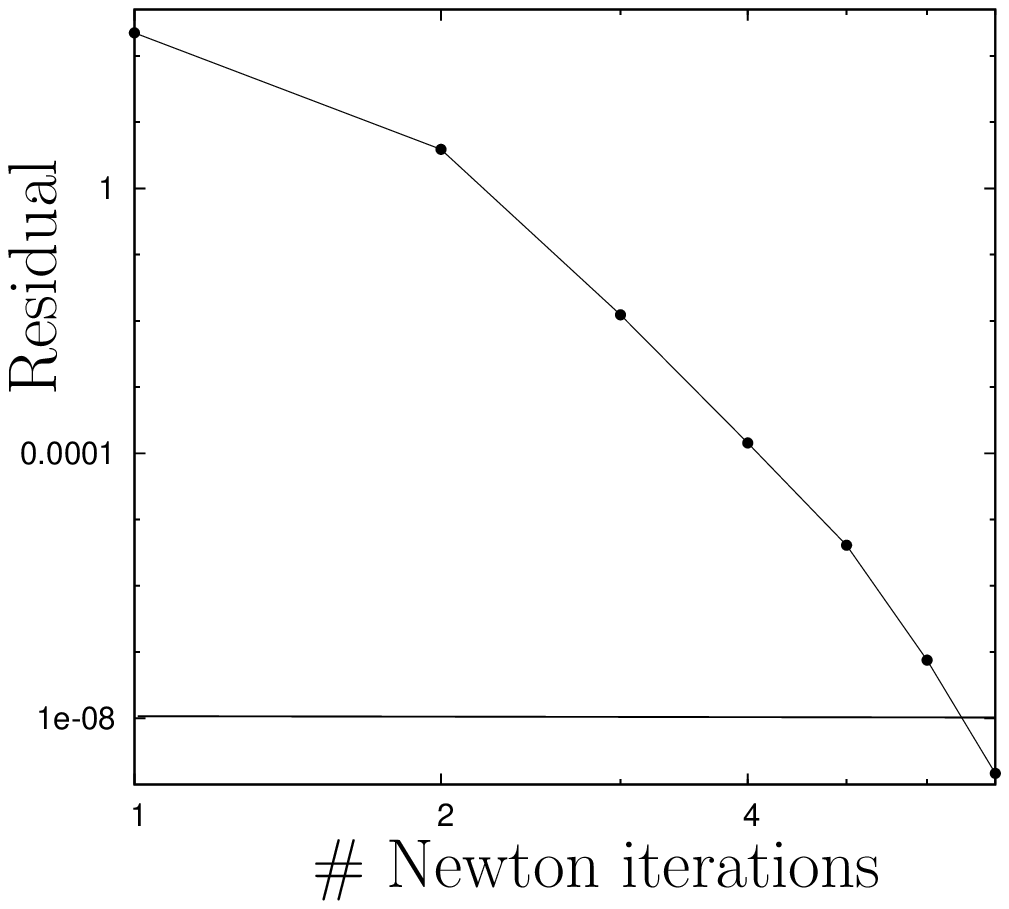}\qquad\includegraphics[width=.45\textwidth]{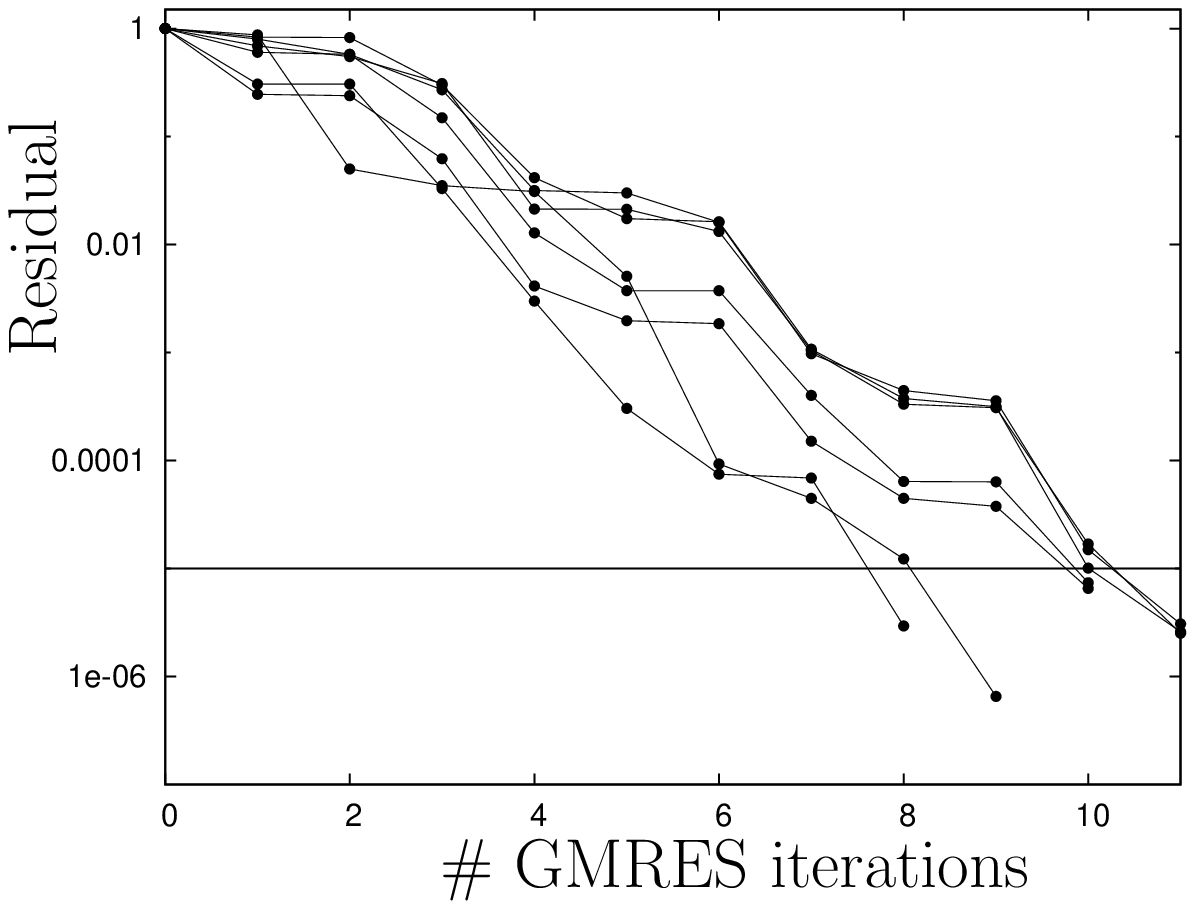}
\end{center}
\caption{\label{fig:NR_GMRES_res}Residuals of the Newton iteration (left) and the corresponding
GMRES iterations (right). The latter is normalised by the norm of the right hand side of Eq.~\ref{eq:periodic_bordered},
i.e. the Newton residual. The tolerance was set at $10^{-8}$ for the Newton iteration and to $10^{-5}$
for the GMRES iteration. Note the super linear convergence of the  former. }
\end{figure}

Starting from a coarse initial approximation, the Newton iterations converged faster than linear,
and each Newton step took between 8 and 11 GMRES iterations, out of a maximum of $N+1=3585$.
Subsequently, we computed the most unstable multipliers, using SLEPc with the {\tt MATSHELL}
that computes products with the bordered matrix, as described in Sec.~\ref{implement:per}.
The most unstable multiplier is $\mu_1=1.111$ and corresponds to a wave number $(1,1)$ perturbation.

\section{Conclusion and future improvements}
In the current paper, we have presented the basic implementation of the model and example computations
to validate it and test its performance. The code will be available publicly \cite{public}.
As it is built on top of PETSc, the user has access to a range of nonlinear and linear solvers and 
preconditioners, which can be used to solve the boundary value problems that typically arise in 
dynamical systems analysis. The periodic orbit computation, presented in Sec.~\ref{ex:per}, is a simple 
example of such a boundary value problem, that has all the ingredients: a module for time-stepping
the system and perturbations and a representation of user-specified, bordered matrices.

The next step in the development of the code is the implementation of pseudo-arclength continuation
of equilibria and periodic orbits. This will enable us, for instance, to complement the bifurcation
diagram of the current test case, Fig.~\ref{fig:bifurcation_diagram}, with the branches of space-dependent
periodic solutions that actually regulate the observed dynamics, in contrast to the highly unstable
spatially homogeneous periodic orbits computed from an ODE reduction of the model.

We expect that our implementation will be useful to researchers 
studying the dynamics of the Liley model, or similar models, such as the model with gap junctions
proposed by Steyn-Ross {\sl et al.} \cite{steynross}. Also, it could be useful for those who
incorporate a similar mean-field model in, for instance, the control of robotic motion
or network models of brain activity.

\section*{Acknowledgements}
LvV was supported by NSERC Grant nr. 355849-2008. Some of the computations were made 
possible by the facilities of the Shared Hierarchical 
Academic Research Computing Network (SHARCNET:www.sharcnet.ca) and Compute/Calcul Canada.

\bibliographystyle{elsarticle-num}
\bibliography{MFM_small}

\end{document}